%% file: KIH_Markov.tex
\documentclass[10pt,a4paper,twoside]{amsart}

\usepackage[utf8]{inputenc}
\usepackage[english]{babel}
\usepackage[T1]{fontenc}
\usepackage{microtype, fancyhdr, lmodern, xcolor}
\usepackage[textsize=9]{todonotes}

\usepackage{amsfonts, amsmath, amsthm, amssymb, mathrsfs, amscd}		
\numberwithin{equation}{section}	
\usepackage{faktor}
\allowdisplaybreaks

\usepackage{url, enumitem}
\usepackage[noadjust]{cite}

\usepackage{multirow,bigdelim, caption}
\usepackage{booktabs}

\usepackage{graphicx}

\usepackage{hyperref}

\theoremstyle{plain}
\newtheorem{thm}{Theorem}[section]
\newtheorem{prop}[thm]{Proposition}

\newtheorem{lem}[thm]{Lemma}

\theoremstyle{remark}
\newtheorem{rmk}[thm]{Remark}

\theoremstyle{definition}
\newtheorem{defn}[thm]{Definition}

\input{Macro}

\begin{document}

\title{Towards a version of Markov's theorem for ribbon torus-links in $\Rr^4$}

\author[C. Damiani]{Celeste Damiani}
\address{Department of Mathematics,
Osaka City University,
Sugimoto, Sumiyoshi-ku,
Osaka 558-8585, Japan}
\email{celeste.damiani@math.cnrs.fr}



\date{\today}
\begin{abstract}
In classical knot theory, Markov's theorem gives a way of describing all braids with isotopic closures as links in~$\mathbb{R}^3$. We present a version of Markov's theorem for extended loop braids with closure in~$B^3 \times S^1$, as a first step towards a Markov's theorem for extended loop braids and ribbon torus-links in~$\Rr^4$. 
\end{abstract}


\maketitle

\section{Introduction}

In the classical theory of braids and links, Alexander's theorem allows us to represent every link as the closure of a braid. Moreover, Markov's theorem states that two braids (possibly with different numbers of strings) have isotopic closures in a $3$-dimensional space if and only if one can be obtained from the other after a finite number of Markov moves, called \emph{conjugation} and \emph{stabilization}. This theorem is a tool to describe all braids with isotopic closures as links in a $3$-dimensional space. Moreover, these two theorems allow us to recover certain link invariants as Markov traces.

When considering \emph{extended loop braid}s as braided annuli in a $4$-dimensional space on one hand, and \emph{ribbon torus-links} on the other hand, we have that a version of Alexander's theorem is a direct consequence of three facts. First of all, every ribbon torus-link can be represented by a welded braid~\cite{Satoh:2000}. Then, a version of Markov's theorem is known for welded braids and welded links~\cite{Kamada:Markov, Kauffman-Lambropoulou:L-move}. Finally, welded braid groups and loop braid groups are isomorphic~\cite{Damiani:Journey}, and loop braids are a particular class of extended loop braids.

In this paper we take a first steps in formulating a version of Markov's theorem for extended loop braids with closure in the space $B^3 \times S^1$. We show that two extended loop braids have closures that are isotopic in $B^3 \times S^1$ if and only if they are conjugate in the extended loop braid groups. The reason for considering \emph{extended} loop braid groups instead of loop braid groups is because this allows to prove a result that is exactly the analogous of the result that we have for $1$-dimensional braids and knots in a $3$-dimensional space. In fact, if we consider two loop braids in the first place, we have that their closures are isotopic as ribbon torus-knots in $B^3 \times S^1$ if and only if the pair of loop braids are conjugate in the extended loop braid group. This is due to the fact that isotopies of ribbon torus-links can introduce a phenomenon called \emph{wen}, which we discuss in Subsection~\ref{SS:extended_ribbon}, on the components of the closed braided objects. Wens are  natural phenomena in the context of ribbon torus-links in~$\Rr^4$, but they are not encoded in the theory of loop braids. Then, extended loop braids, who encode wens, seem to be the most natural analogue of classical braids, and the most appropriate notion that we need to consider.

\subsection{Structure of the paper}
In Section~\ref{S:ExtLoopBraids} we give an overview of the many equivalent interpretations of extended loop braid groups, which are the braided objects coming to play in our main result. A particular focus will be placed on the definition of extended loop braids as braided annuli in a $4$-dimensional space. When we want to make clear that we are using this interpretation for extended loop braids, we use the terminology \emph{ribbon braids}. We recall several results on these objects, and we use them to prove that every ribbon braid can be parametrized by a normal isotopy (Proposition~\ref{L:NormalParametrization}).

In Section~\ref{S:RibbonLinks} we introduce the knotted counterpart of ribbon braids, which are ribbon torus-links. 

In Section~\ref{S:Markov} we present the main result of this paper. This is a version of Markov's theorem for ribbon torus-links living in the space $B^3 \times S^1$ (Theorem~\ref{T:toro}). 

Finally, in Section~\ref{S:future} we discuss possible ideas to complete the main result of this paper to a complete Markov's theorem for ribbon torus-links in~$\Rr^4$.

\section{Extended loop braid groups and their equivalent definitions}
\label{S:ExtLoopBraids}

\emph{Loop braid groups} were introduced under this name for the first time by Xiao-Song Lin in 2007~\cite{Lin:2008}, although they had been considered before in other contexts and with other terminologies, for instance \emph{groups of basis-conjugating automorphisms} in \cite{Savushkina:1996} and \emph{welded braid groups} in~\cite{Fenn-Rimanyi-Rourke:1997}. 

The groups we call \emph{extended loop braid groups} appeared sooner in the literature, in~\cite{Dahm} and~\cite{Goldsmith:MotionGroups}, who called them \emph{motion groups of a trivial link of unknotted circles  in $\Rr^3$}, but since then have been less treated in the literature.

In terms of configuration spaces, both groups appear in~\cite{BrendleHatcher:2013}, loop braid groups as \emph{untwisted ring groups}, and extended loop braid groups as~\emph{ring groups}. In this paper we focus on extended loop braid groups, for which we choose to adapt Lin's notation because it gives a good visual idea of the considered objects, while being more compact. In fact, the elements of both these groups can be seen as trajectories travelled by loops as they move in a $3$-dimensional space to exchange their positions under some admissible motions. The ``extended'' attribute highlights the fact that in extended loop braid groups we admit an extra motion that can be described  as a $180^\circ$-flip of a loop.
For a detailed survey on loop braid groups, extended loop braid groups and the explicit equivalences among the different definitions, we refer to~\cite{Damiani:Journey}.

We dedicate this section to recall several definitions of extended loop braid groups, and give the terminology used in the different contexts. The diversity of points of view will be useful in the proof of the main result of this paper (Theorem~\ref{T:toro}), since it provides many approaches and tools to tackle problems involving extended loop braid groups and other knotted objects in the $4$-dimensional space.

\subsection{Extended loop braids as mapping classes}
\label{SS:Mapping classes}
We present here a first definition for extended loop braid groups in terms of mapping classes of a $3$-ball with $n$ circles that are left setwise invariant in its interior.

Let us fix $\nn \in \Nn$, and let $C = C_1 \sqcup \cdots \sqcup C_\nn$ be a collection of $n$ disjoint, unknotted, oriented circles, that form a trivial link of $\nn$ components in the interior of the $3$-ball~$B^3$.  A self-homeomorphism of the pair $(B^3,C)$ is an homeomorphism $f \colon B^3 \to B^3$ that fixes $\partial B^3$ pointwise, preserves orientation on~$B^3$, and globally fixes~$C$. Every self-homeomorphism of $(B^3, C)$ induces a permutation on the connected components of~$C$ in the natural way.
We consider the mapping class group of $B^3$ with respect to~$C$ to be the group of isotopy classes of self-homeomorphisms of~$(B^3, C)$, with multiplication determined by composition. We denote it by~$\MCG {B^3} {C^\ast}$. 

\begin{rmk}
The ``$^\ast$'' on the submanifold~$C$ is to indicate that homeomorphisms do not preserve the orientation of the connected components of~$C$. This is the difference between extended loop braid groups and loop braid groups in this context. In fact, in the latter, the  homeomorphisms preserve orientation on~$C$.
\end{rmk}

\begin{rmk}
\label{R:pi0}
A map $f$ from a topological space $X$ to $\Homeo(B^3; C^\ast)$ is continuous if and only if the map $X \times B^3 \to B^3$ sending $(x, y) \mapsto f(x)(y)$ is continuous~\cite{Kelley:1975}. Taking $X$ equal to the unit interval~$I$, we have that two self-homeomorphisms are isotopic if and only if they are connected by a path in~$\Homeo(B^3; C^\ast)$. Therefore~$\MCG {B^3} {C^\ast} =  \pi_0(\Homeo(B^3; C^\ast))$. The same can be said for the pure groups, $\PMCG  {B^3} {C^\ast}= \pi_0(\PHomeo(B^3; C^\ast))$. 
\end{rmk}

\begin{defn}
\label{D:LBN}
For $\nn \geq 1$, the \emph{extended loop braid group}, denoted  by $\LBE\nn$, is the mapping class group~$\MCG{B^3}{C^\ast}$. 
\end{defn}

\subsection{Extended loop braids as loops in a configuration space}
\label{SS:Configurations}
The second interpretation of extended loop braid groups $\LBE\nn$ that we give is in terms of configuration spaces, and has been introduced in~\cite{BrendleHatcher:2013}. Let $\nn \geq 1$, and consider the space of configurations of $\nn$ Euclidean, unordered, disjoint, unlinked circles in $B^3$, denoted by~$\Ri\nn$. The \emph{ring group} $\R\nn$ is its fundamental group.
 
Remark that in Subsection~\ref{SS:Mapping classes} we were not considering Euclidean circles as moving objects, but the components of a trivial link. We shall see now that these two families of objects are deeply related. Let~$\mathcal{L}_{\nn}$ be the space of configurations of smooth trivial links with $\nn$ components in $\Rr^3$: the following result allows us to consider the fundamental group of~$\mathcal{L}_\nn$ as being isomorphic to~$\R\nn$. 

\begin{thm}[{\cite[Theorem~1]{BrendleHatcher:2013}}]
\label{T:RelaxingCircles}
For $\nn \geq 1$, the inclusion of $\Ri\nn$ into~$\mathcal{L}_{\nn}$ is a homotopy equivalence.
\end{thm}

As anticipated, the groups $\R\nn$ are isomorphic to the groups~$\LBE\nn$, as stated in the next theorem. Its proof heavily relies on Wattenberg's results~\cite[Lemma 1.4 and Lemma 2.4]{Wattenberg:1972} implying that the topological mapping class groups of the $3$-ball with respect to an $\nn$-components trivial link are isomorphic to the $C^\infty$-mapping class groups of the same pair. In other terms, we have that $\pi_0(\Homeo(B^3; C^\ast)) \cong \pi_0(\Diffeo(B^3; C^\ast))$. We can define an evaluation map from $\Diffeo(B^3)$ to the space of configurations of a smooth trivial link with $\nn$ \emph{ordered} components in $\Rr^3$, that we denote by~$\mathcal{PL}_{\nn}$.  We can refer to $\mathcal{PL}_{\nn}$ as to the \emph{pure} configuration space of a smooth trivial link. Fixed the $\nn$ components of a trivial link in the interior of the $3$-ball, this evaluation map sends self-diffeomorphisms of $B^3$ to the image of the $\nn$ components through the considered self-diffeomorphism: 
\begin{equation}
\label{E:evaluation}
\ep \colon \Diffeo(B^3) \longrightarrow \mathcal{PL}_{\nn}.
\end{equation}
This map can be proved to be a locally trivial fibration~\cite[Lemma~3.8]{Damiani:Journey}. This fibration is then used as the main ingredient to prove the following, through the construction of exact sequences and a commutative diagram.

\begin{thm}[{\cite[Theorem~3.10]{Damiani:Journey}}]
\label{T:PureRing}
For $\nn \geq 1$, there is a natural isomorphism between ring group $\R\nn$ and the extended loop braid group~$\LBE\nn$.
\end{thm}

Brendle and Hatcher, in \cite[Proposition~3.7]{BrendleHatcher:2013}, give a presentation for the ring groups~$\R\nn$, and so, for~$\LBE\nn$.

\begin{prop}
\label{P:PresLBE}
For $\nn \geq 1$, the group $\LBE\nn$ admits the presentation given by generators $\{\sig\ii, \rr\ii \mid \ii=1, \dots , \nno \}$ and $\{\tau_\ii \mid \ii=1, \dots , \nn \}$, subject to relations:
\begin{equation}
\label{E:Rpresentation}
\begin{cases}
\sig{i} \sig j = \sig j \sig{i}  \, &\text{for } \vert  i-j\vert > 1\\
\sig{i} \sig {i+1} \sig{i} = \sig{i+1} \sig{i} \sig{i+1} \, &\text{for } i=1, \dots, \nn-2 \\
\rr{i} \rr j = \rr j \rr{i}  \, &\text{for }  \vert  i-j\vert > 1\\
\rr\ii\rr{i+1}\rr\ii = \rr{i+1}\rr\ii\rr{i+1} \, &\text{for }  i=1, \dots, \nn-2  \\
\rrq{i}2 =1 \, &\text{for }  i=1, \dots, \nno \\
\rr{i} \sig{j} = \sig{j} \rr{i}   \, &\text{for }  \vert  i-j\vert > 1\\
\rr{i+1} \rr{i} \sig{i+1} = \sig{i} \rr{i+1} \rr{i} \,  &\text{for }  i=1, \dots, \nn-2  \\
\sig{i+1} \sig{i} \rr{i+1} = \rr{i} \sig{i+1} \sig{i}  \,  &\text{for }  i=1, \dots, \nn-2 \\
\tau_{i} \tau_j = \tau_j \tau_{i}  \, &\text{for }    i \neq j \\
\tau_\ii^2=1 \, &\text{for }  i=1, \dots, \nn \\
\sig{i} \tau_j = \tau_j \sig{i}  \, &\text{for }  \vert  i-j\vert > 1 \\
\rr{i} \tau_j = \tau_j \rr{i}  \, &\text{for }  \vert  i-j\vert > 1 \\
\tau_i \rr\ii = \rr\ii \tau_{i+1} \, &\text{for }  i=1, \dots, \nno  \\
\tau_\ii \sig\ii = \sig\ii \tau_{i+1}  \, &\text{for } i=1, \dots, \nno \\
\tau_{i+1} \sig\ii = \rr\ii \siginv\ii \rr\ii \tau_\ii  \, &\text{for } i=1, \dots, \nno. \\
\end{cases}
\end{equation}
\end{prop}

\begin{figure}[hbtp]
\centering
\includegraphics[scale=.5]{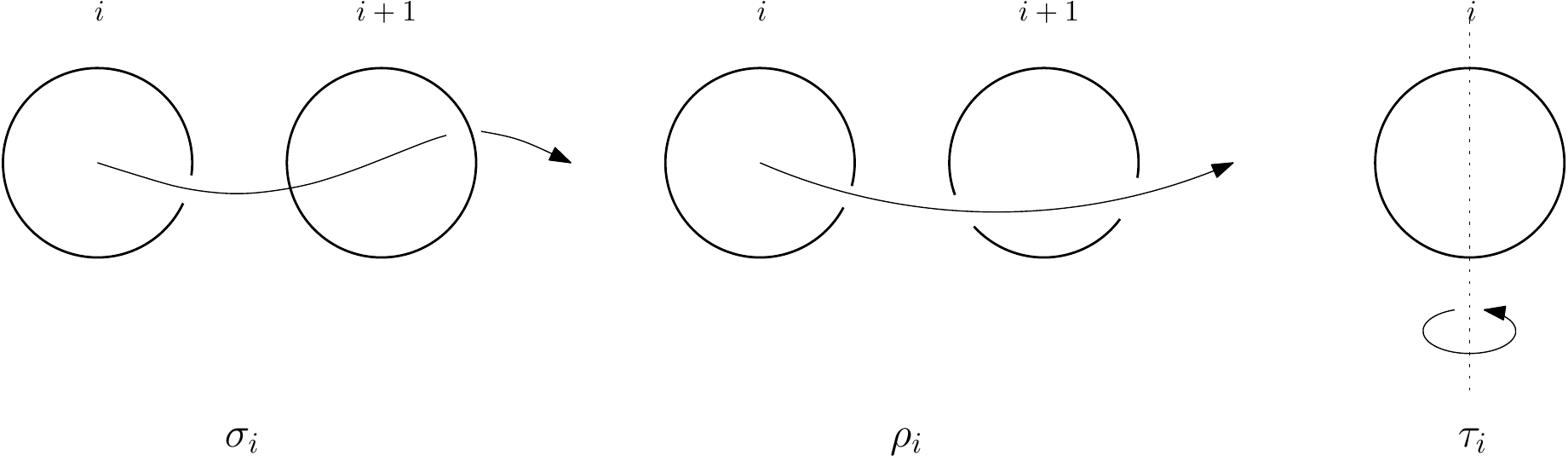}
\caption{Elements $\sig\ii$, $\rr\ii$ and $\tau_\ii$.}
\label{F:Flips}
\end{figure}

The elements~$\sig\ii$,~$\rr\ii$, and $\tau_\ii$ of the presentation represent the following loops in $\R\nn$: if we place the $\nn$ rings in a standard position in the $yz$-plane with centers along the $y$-axis, then $\sig\ii$ is the loop that permutes the $\ii$-th and the $(\ii+1)$-st circles by passing the $\ii$-th circle through the $(\ii+1)$-st; $\rr\ii$ permutes them passing the $\ii$-th around the $(\ii+1)$-st, and $\tau_i$ is the loop that flips by $180^\circ$ the $i$-th circle, see Figure~\ref{F:Flips}.

\subsection{Extended loop braids as automorphisms of the free groups}
We now give an interpretation of extended loop braids in terms of automorphisms of~$\F\nn$, the free groups of rank~$\nn$. 
Fixing $\nn \geq 1$, we consider the  automorphisms that send each generator of $\F\nn$ to a conjugate of some generator, or its inverse: these are in bijection with elements of~$\LBE\nn$. We start recalling a result of Dahm's unpublished thesis~\cite{Dahm}, that appears in the last section of Goldsmith's paper~\cite{Goldsmith:MotionGroups}. 

\begin{thm}[{\cite[Theorem~5.3]{Goldsmith:MotionGroups}}]
\label{T:Gold}
For $\nn \geq 1$, there is an injective map from the extended loop braid group $\LBE\nn$ into~$\Aut(\F\nn)$, where $\F\nn$ is the free group on $\nn$ generators $\{x_1, \dots, x_\nn\}$, and its image is the subgroup~$\PC\nn$, consisting of all automorphisms of the form $\alpha \colon x_\ii \mapsto w_\ii\inv x^{\pm 1}_{\pi({\ii})} w_\ii$ where $\pi$ is a permutation and $w_\ii$ is a word in~$\F\nn$. Moreover, the group $\PC\nn$ is generated by the automorphisms $\{ \sig1, \dots \sig\nno,\rr1, \dots \rr\nno, \tau_1, \dots, \tau_\nn\}$ defined as: 
\begin{align}
\label{E:sigma}
\sig\ii &: \begin{cases}
            x_\ii \mapsto x_{\ii+1}; &\\
            x_{\ii+1} \mapsto x_{\ii+1}\inv x_\ii x_{\ii+1}; &\\
            x_\jj \mapsto x_\jj, \ &\text{for} \ \jj \neq \ii, \ii+1.
        \end{cases} \\
\label{E:rho}
\rr\ii &: \begin{cases}
            x_\ii \mapsto x_{\ii+1}; & \\
            x_{\ii+1} \mapsto x_\ii; &\\
            x_\jj \mapsto x_\jj, \ &\text{for} \ \jj \neq \ii, \ii+1.
		\end{cases} \\
\label{E:tau}
\tau_\ii &: \begin{cases}
            x_\ii \mapsto x\inv_\ii; &\\
            x_\jj \mapsto x_\jj, \ &\text{for} \ \jj \neq \ii.
        \end{cases}
\end{align}
\end{thm}

This result is the analogue of Artin's characterization of usual braids as automorphisms of the free group. In an intuitive way, we use for the automorphisms of $\PC\nn$ the notations of the corresponding elements of the mapping class group\footnote{In~\cite{Damiani:Journey} these groups are denote by~$\PC\nn^\ast$, while $\PC\nn$ is used for the groups of automorphisms of the form~$\alpha \colon x_\ii \mapsto w_\ii\inv x_{\pi({\ii})} w_\ii$.}.

In~\cite{Fenn-Rimanyi-Rourke:1997} Fenn, Rim\'{a}nyi and Rourke consider the subgroups of $\Aut(\F\nn)$ generated only by the sets of elements $\{\sig\ii \mid i=1, \dots \nn-1\}$ and $\{\rr\ii \mid i=1, \dots \nn-1\}$. They call these groups by the name \emph{braid-permutation groups}, and they prove independently from Dahm and Goldsmith that they are isomorphic to the groups  of all automorphisms of $\Aut(\F\nn)$ of the form $\alpha \colon x_\ii \mapsto w_\ii\inv x_{\pi({\ii})} w_\ii$ where $\pi$ is a permutation and $w_\ii$ is a word in~$\F\nn$.

\subsection{Extended loop braids as ribbon braids}
\label{SS:extended_ribbon}
The next interpretation of extended loop braids will be the one that we will focus on in the main result of this paper. This is an approach in terms of braided objects in a $4$-dimensional space. Extended loop braids in this context are called \emph{ribbon braids}, when we want to specify the used interpretation\footnote{In the survey~\cite{Damiani:Journey} the terminology \emph{ribbon braids} refers to loop braids seen as braided objects in the $4$-dimensional braid, while the terminology \emph{extended ribbon braids} refers to extended loop braids. We chose to simplify.}.

We need some notation before giving the definition of ribbon braids and their equivalence to extended loop braids. Let $\nn \geq 1$, and let $D_1, \dots, D_\nn$ be a collection of disks in the $2$-ball~$B^2$. Let $C_i =\partial D_i$ be the oriented boundary of~$D_i$.
Let us consider the $4$-ball $B^4 \cong B^3 \times I$, where $I$ is the unit interval. For any submanifold $X \subset B^m \cong B^{m-1} \times I$, with~$m=3,4$, we use the following dictionary. To keep the notation readable, here we denote the interior of a topological space by ``$\mathrm{int}(\phantom{x} )$'', whereas anywhere else it is denoted by ``$\mathring{\phantom{x}}$''.
\begin{itemize}
\item $\partial_\ep X = X \cap (B^{m-1} \times \{\ep\})$, with~$\ep \in \{0, 1\}$;
\item $\partial_\ast X = \partial X \setminus \Big( \mathrm{int}(\partial_0 X) \sqcup \mathrm{int}(\partial_1 X) \Big)$;
\item $\overset{*}{X}=X \setminus \partial_\ast X$.
\end{itemize}

The image of an immersion $Y \subset X$ is said to be \emph{locally flat} if and only if it is locally homeomorphic to a linear subspace $\Rr^k$ in $\Rr^m$ for some $k \leq m$, except on $\partial X$ and/or $\partial Y$, where one of the $\Rr$ summands should be replaced by $\Rr_+$. Let $Y_1, Y_2$ be two submanifolds of $B^m$. The intersection $Y_1 \cap Y_2 \subset X$ is called \emph{flatly transverse} if and only if it is locally homeomorphic to the transverse intersection of two linear subspaces $\Rr^{k_1}$ and $\Rr^{k_2}$ in $\Rr^m$ for some positive integers $k_1, k_2 \leq m$ except on $\partial X$, $\partial Y_1$ and/or $\partial Y_2$, where one of the $\Rr$ summands should be replaced by $\Rr_+$. In the next definition we introduce the kind of singularities we consider.

\begin{defn}
\label{D:RibbonDisk}
Let $Y_1, Y_2$ be two submanifolds of $B^4$. \emph{Ribbon disks} are intersections $D = Y_1 \cap Y_2$ that are isomorphic to the $2$-dimensional disk,  such that $D \subset \mathring{Y_1}$, $\mathring{D} \subset\mathring{Y_2}$ and $\partial D $ is an essential curve in~$\partial Y_2$.
\end{defn}
These singularities are the $4$-dimensional analogues of the classical notion of ribbon singularities introduces by Fox in~\cite{Fox:Ribbon}.

\begin{defn}
\label{D:ribbonBraid}
Let $A_1, \dots, A_\nn$ be locally flat embeddings in~$\stackrel{*}{B^4}$ of $\nn$ disjoint copies of the oriented annulus $S^1 \times I$. We say that   
\[b= \bigsqcup_{i\in \{1, \dots, n \}} A_i \]
is a \emph{geometric ribbon braid} if:
\begin{enumerate}
\item \label{bordo} the boundary of each annulus $\partial A_i$ is a disjoint union  $C_i \sqcup C_j$, for $C_i \in \partial_0 B^4$ and for some $C_j \in \partial_1 B^4$. The orientation induced by $A_i$ on $\partial A_i$ coincides with the one of the two boundary circles $C_i$ and~$C_j$;
\item the annuli $A_i$ are fillable, in the sense that they bound immersed $3$-balls $\subset \Rr^4$ whose singular points consist in a finite number of ribbon disks;
\item \label{trans} it is transverse to the lamination $\bigcup_{t \in I} B^3 \times \{t\}$ of~$B^4$, that is: at each parameter~$t$, the intersection between $b$ and $B^3 \times{t}$ is a collection of exactly $\nn$ circles;
\end{enumerate}

The group of \emph{ribbon braids}, denoted by $\rB\nn$, is the group of equivalence classes of geometric ribbon braids up to continuous deformations through the class of geometric ribbon braids fixing the boundary circles, equipped with the natural product given by stacking and reparametrizing. The unit element for this product is the \emph{trivial ribbon braid}~$U= \bigsqcup_{i\in \{1, \dots, n \}} C_i \times [0, 1]$.
\end{defn}

The monotony condition allows us to consider the interval $I$ in $B^4 = B^3 \times I$ as a time parameter, and to think of a ribbon braid as a trajectory $\beta=\big( C_1(t), \dots, C_\nn(t)\big)$ of circles in $B^3 \times I$. This trajectory corresponds to a parametrization of the ribbon braid. This interpretation is also referred to in terms of \emph{flying rings} in~\cite{BarNatan-Dancso:Survey}.  When one of the $\nn$ circles that we have at each time $t$ makes a half-turn, we have what is called a \emph{wen} on the corresponding component. One can think of a wen as an embedding in $\Rr^4$ of a Klein bottle cut along a meridional circle. A detailed treatment of wens can be found in Kanenobu and Shima's paper~\cite{KanenobuShima:2002}.

The following result states the equivalence of the interpretations of $\LBE\nn$ as mapping class groups and as ribbon braid groups. Its proof consists in explicitly defining an isomorphism between $\rB\nn$ and $\R\nn$, and composing it with the isomorphism from Theorem~\ref{T:PureRing}.

\begin{thm}[{\cite[Theorem~5.17]{Damiani:Journey}}]
\label{T:ribbon-e-loop}
For $\nn \geq 1$, there is an isomorphism between the ribbon braid group $\rB\nn$ and the extended loop braid group~$\LBE\nn$ .
\end{thm}

We can show that when two ribbon braids are equivalent in the sense of Definition~\ref{D:ribbonBraid}, there is an ambient isotopy of $\Rr^4$ bringing one to the other.
\begin{thm}[{\cite[Theorem~5.5]{Damiani:Journey}}]
Every relative isotopy of a geometric ribbon braid in $B^3 \times I$ extends to an isotopy of $B^3 \times I$ in itself constant on the boundary.
\end{thm}

This result is true also for \emph{surface links}, which are closed surfaces locally flatly embedded in~$\Rr^4$~\cite[Theorem~6.7]{Kamada:libro}.

With the results we recalled, we prove now that given a geometric ribbon braid $b$ and its set of starting set of circles, we can find a normal isotopy parametrizing it.
We take $C = (C_{1}, \dots, C_{n})$ to be an ordered tuple of $\nn$ disjoint, unlinked, unknotted circles living in~$B^3$.
We consider the space of configurations of ordered smooth trivial links of $\nn$ components $\mathcal{PL}_\nn$ introduced in Subsection~\ref{SS:Configurations}.
As mentioned above, we have an evaluation map 
\begin{equation*}
\ep \colon \Diffeo(B^3) \longrightarrow \mathcal{PL}_{\nn}
\end{equation*}
sending a self-diffeomorphism $f$ to~$f(C)$. We remark that $f(C)$ is an ordered tuple of $\nn$ disjoint, unlinked, trivial, smooth knots living in~$B^3$, which is a locally trivial fibration with fibre the group of self-diffeomorphisms of the pair $(B^3, C)$ that send each connected component of $C$ to itself.
Composing $\ep$ with the covering map $\mathcal{PL}_{\nn} \to \mathcal{L}_{\nn}$, seeing $\mathcal{L}_{\nn}$ as the orbit space with of the action of the symmetric group of $\mathcal{PL}_{\nn}$, we define a locally trivial fibration
\[
\tilde{\ep}\colon \Diffeo(D^3) \longrightarrow \mathcal{L}_{\nn}
\]
sending $f$ to~$f(C)$. More details on this construction can be found in~\cite{Damiani:Journey}.


\begin{lem}
\label{L:NormalParametrization}
Let $\nn\geq 1$. For every geometric ribbon braid $b \subset B^4$ on $n$ components, there is a normal isotopy parametrizing $b$.
\end{lem}
\proof

Let us consider a geometric ribbon braid~$b$, through the isomorphism between $\rB\nn$ and $\R\nn$ (Theorem~\ref{T:ribbon-e-loop}). This gives rise to a loop $f^b \colon  I \to \mathcal{PL}_{\nn} \subset \mathcal{L}_{\nn}$ sending $t \in I$ into the unique $n$-circles set $b_t$ such that 
\[
b\cap(B^3 \times {I})=b_t \times \{t\}.
\]
This loop begins and ends at the point $\tilde{\ep}(\id_{B^3}) \in \mathcal{L}_{\nn}$ represented by~$C$. Being $\tilde{\ep}$ a fibration, we apply the homotopy lifting property, and lift $f^b$ to a path $\hat{f^b}\colon I \to \Diffeo(B^3)$ beginning at $\tilde{\ep}\inv(C)=\Diffeo(B^3; C^\ast)$ and ending at~$\id_{B^3}$. The path $\hat{f^b}$ is a normal isotopy. The commutativity $\tilde{\ep} \circ \hat{f^b}=f^b$ means that this isotopy parametrizes~$b$.
\endproof

\subsection{Extended loop braids as extended welded braids}
\label{SS:ext_welded_braids}
\begin{figure}[hbt]
\centering
\includegraphics[scale=0.6]{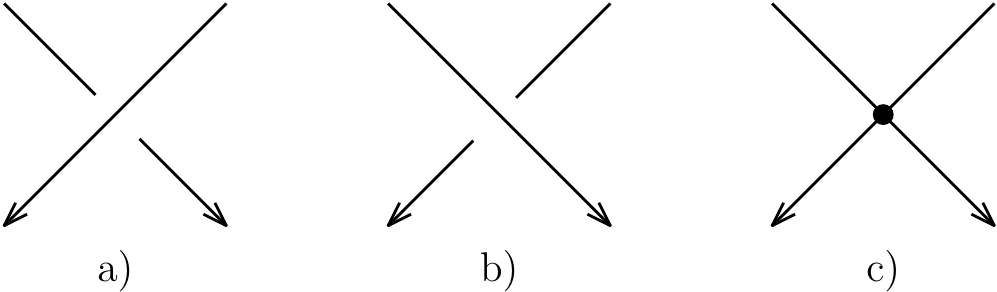}
\caption{a) Classical positive crossing, b) Classical negative crossing, c) Welded crossing.}
\label{F:Crossings}
\end{figure}

\begin{figure}[htb]
\centering
\includegraphics[scale=0.6]{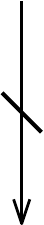}
\caption{A wen mark on a strand.}
\label{F:WenMark}
\end{figure}

In this part we discuss $1$-dimensional diagrams immersed in a $2$-dimensional space for extended loop braids. An \emph{extended welded braid diagram} on $\nn$ strings is a planar diagram composed by a set of $\nn$ oriented and monotone $1$-manifolds immersed in $\Rr^2$ starting from $\nn$ points on a horizontal line at the top of the diagram down to a similar set of $\nn$ points at the bottom of the diagram. The $1$-manifolds are allowed to cross in transverse double points, which will be decorated in three kinds of ways, as shown in Figure~\ref{F:Crossings}. Depending on the decoration, double points will be called: \emph{classical positive} crossings, \emph{classical negative} crossings and \emph{welded} crossings. On each $1$-manifold there can possibly be marks as in~Figure~\ref{F:WenMark}, which we will call \emph{wen marks}.

\begin{figure}[hbt]
\centering
\includegraphics[scale=0.6]{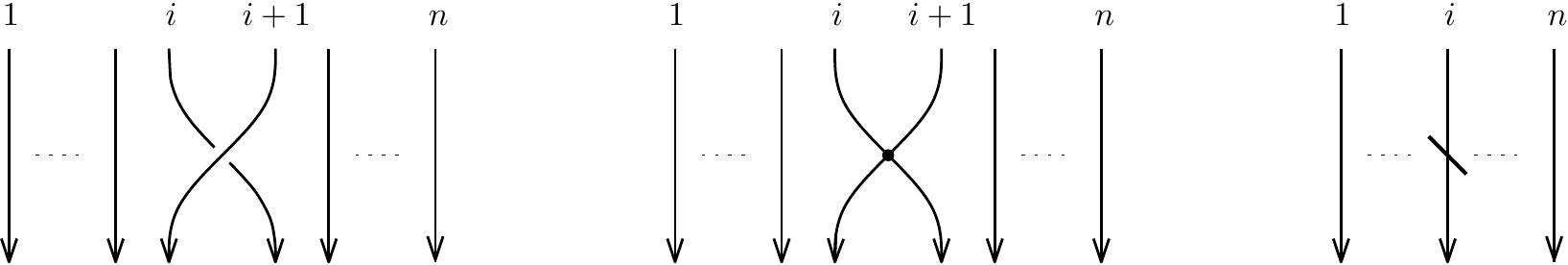}
\caption{Elementary diagrams $\sig\ii$, $\rr\ii$, and~$\tau_\ii$.}
\label{F:SigRhoTau}
\end{figure}

Let us assume that the double points occur at different $y$-coordinates. Then an extended welded braid diagram determines a word in the elementary diagrams illustrated in Figure~\ref{F:SigRhoTau}. We call $\sig\ii$ the elementary diagram representing the $(i+1)$-th strand passing over the $i$-th strand, $\rr\ii$ the welded crossing of the strands $i$ and~$(i+1)$, and $\tau_\ii$ the wen mark diagram.

\begin{figure}[hbt]
\centering
\includegraphics[scale=.6]{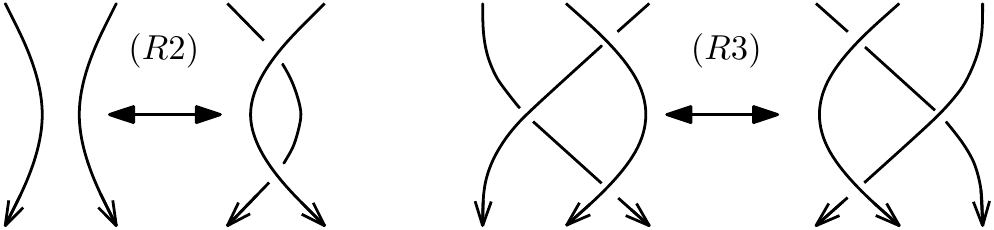}
\caption{Classical Reidemeister moves for braid-like objects.}
\label{F:Classical}
\end{figure}

\begin{figure}[hbt]
\centering
\includegraphics[scale=.6]{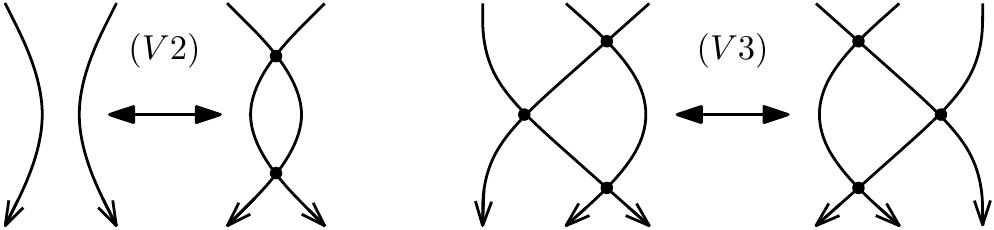}
\caption{Virtual Reidemeister moves for braid-like objects.}
\label{F:Virtual}
\end{figure}

\begin{figure}[hbt]
\centering
\includegraphics[scale=.6]{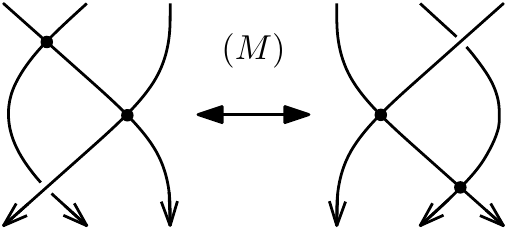}
\caption{Mixed Reidemeister moves.}
\label{F:Mixed}
\end{figure}

\begin{figure}[hbt]
\centering
\includegraphics[scale=.6]{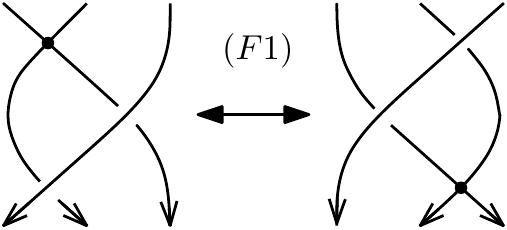}
\caption{Welded Reidemeister moves.}
\label{F:Welded}
\end{figure}

\begin{figure}[htb]
\centering
\includegraphics[scale=0.6]{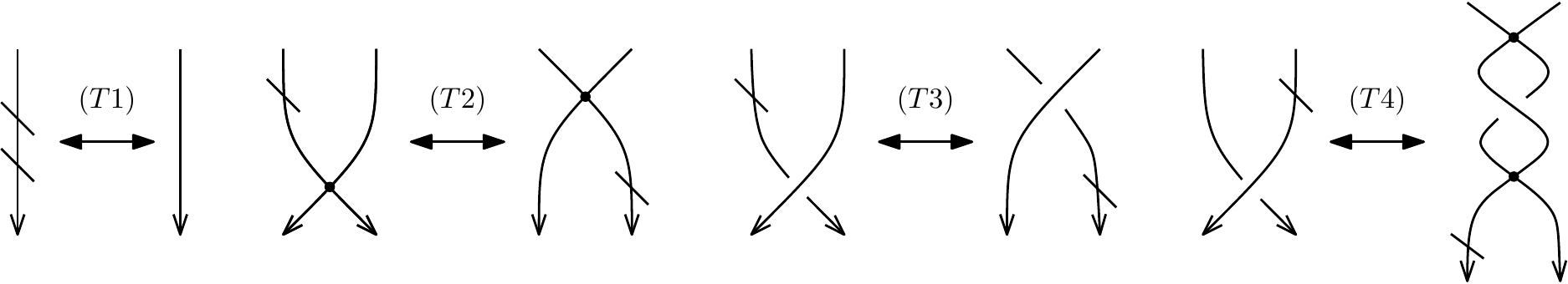}
\caption{Extended Reidemeister moves.}
\label{F:ReidWen}
\end{figure}

\begin{defn}
\label{D:bWelded}
An \emph{extended welded braid} is an equivalence class of extended welded braid diagrams under the equivalence relation given by isotopy of $\Rr^2$ and the following moves: 
\begin{itemize}
\item classical Reidemester moves (Figure~\ref{F:Classical});
\item virtual Reidemeister moves (Figure~\ref{F:Virtual});
\item mixed Reidemeister moves (Figure~\ref{F:Mixed});
\item welded Reidemeister moves (Figure~\ref{F:Welded});
\item extended Reidemester moves (Figure~\ref{F:ReidWen}).
\end{itemize} 
This equivalence relation is called \emph{(braid) generalized Reidemeister equivalence}. 
For~$\nn \geq 1$, the \emph{extended welded braid group} on $\nn$ strands $\WBE\nn$ is the group of equivalence classes of extended welded braid diagrams by generalized Reidemeister equivalence. The group structure on these objects is given by: stacking and rescaling as product, braid mirror image as inverse, and the trivial diagram as identity.
\end{defn}

\begin{rmk}
If wen marks were not allowed, the group defined would be the group of \emph{welded braids}~$\WB\nn$, introduced by Fenn, Rim\'{a}nyi and Rourke in~\cite{Fenn-Rimanyi-Rourke:1997}. This group is isomorphic to \emph{loop braid groups}~$\LB\nn$.
\end{rmk}

In~\cite{Satoh:2000} the author defines a surjective map $Tube$ from welded knotted objects to ribbon knotted objects in dimension~$4$. 
This map $Tube$ can be easily defined on extended welded braids and and extended loop braids, in their interpretation as ribbon braids. Its definition uses as a stepping stone a projection of ribbon braids onto a certain class of $3$-dimensional surfaces, called \emph{broken surface diagrams}. We do not treat them in this paper since they are not relevant to the main result. However they are an interesting way of representing ribbon braids, and more detail can be found in~\cite{BaezWiseCrans:Exotic}. In the framework of extended welded braids and ribbon braids, it can be proved that the $Tube$ map is an isomorphism~\cite[Theorem~6.12]{Damiani:Journey}.  Hence, we have the last isomorphism that we recall in this overview on extended loop braids. 
\begin{thm}
For $\nn \geq 1$, there is an isomorphism between the extended welded braid group $\WBE\nn$ and the extended loop braid group~$\LBE\nn$.
\end{thm}

\subsection{Pure subgroups}
As in the case of classical braid groups $\BB\nn$, we have a notion of pure subgroups for the extended loop braid groups $\LBE\nn$. Let us consider the first definition we gave for extended loop braids, as elements of $\MCG {B^3}{C^\ast}$, where $C = C_1 \sqcup \cdots \sqcup C_\nn$ is a collection of $n$ disjoint, unknotted, oriented circles, that form a trivial link of $\nn$ components. Let $p\colon \LBE\nn \to S_\nn$ be the homomorphism that forgets the details of the braiding, remembering only the permutation of the circles. Then the \emph{pure extended loop braid group} $\PLBE\nn$ is the kernel of~$p$. In each one of the approaches to extended loop braid groups that we exposed, such subgroups can be defined with tools inherent to the particular context. 
We will not dwell on these groups here, but they are discussed in all the references we gave on extended loop braid groups throughout this section.

\section{Ribbon torus-links}

\label{S:RibbonLinks}
In this part we introduce the knotted counterpart of extended loop braid groups: ribbon torus-links. Classical references for these objects are~\cite{Kamada:libro, Kawauchi:Survey, Yajima:2manifolds}.

\begin{defn}
\label{D:ribbonknot}
A \emph{geometric ribbon torus-knot} is an embedded oriented torus $S^1 \times S^1 \subset \Rr^4$ which is \emph{fillable}, in the sense that it bounds a \emph{ribbon torus}, \ie, an oriented immersed solid torus $D^2 \times S^1\subset \Rr^4$ whose singular points consist in a finite number of ribbon disks. \emph{Ribbon torus-knots} are equivalence classes of geometric ribbon torus-knots defined up to ambient isotopy.
\end{defn}

\begin{rmk}
\label{R:wen}
Wens can appear on portions of a ribbon knot, but for an argument of coherence of the co-orientation, there are an even number of them on each component, and they cancel pairwise, as remarked in~\cite[proof of Proposition 2.4]{Audoux:2016}.
\end{rmk}

\begin{defn}
\label{D:ribbolink}
A \emph{geometric ribbon torus-link} with $\nn$ components is the embedding of a disjoint union of $\nn$ oriented fillable tori. The set of \emph{ribbon torus-links} is the set of equivalence classes of geometric ribbon torus-knots defined up to ambient isotopy.
\end{defn}

\subsection{Extended welded diagrams for ribbon torus-links}

An \emph{extended welded link diagram} is the immersion in~$\Rr^2$ of a collection of disjoint, closed, oriented $1$-manifolds such that all multiple points are transverse double points. Double points are decorated with classical positive, classical negative, or welded information as in Figure~\ref{F:Crossings}. On each $1$-manifod there can possibly be an even number of wen marks as in~Figure~\ref{F:WenMark}, the motivation for this lying in Remark~\ref{R:wen}. We assume that extended welded link diagrams are the same if they are isotopic in~$\Rr^2$. 
Taken an extended welded link diagram $K$, we call \emph{real crossings} its set of classical positive and classical negative crossings.

\begin{figure}[htb]
	\centering
		\includegraphics[scale=0.7]{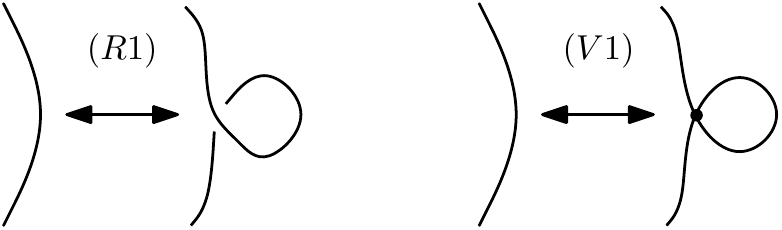} 
	\caption{Reidemeister moves of type I.}
	\label{F:ReidOne}
\end{figure}

\begin{defn}
\label{D:ExWeldedLink}
An \emph{extended welded link} is an equivalence class of extended welded link diagrams under the equivalence relation given by isotopies of $\Rr^2$, moves from Definition~\ref{D:bWelded}, and classical and virtual Reidemeister moves $(R1)$ and $(V1)$ as in Figure~\ref{F:ReidOne}. This equivalence relation is called \emph{generalized Reidemeister equivalence}. 
\end{defn}

\begin{figure}[htb]
\centering
\includegraphics[scale=0.6]{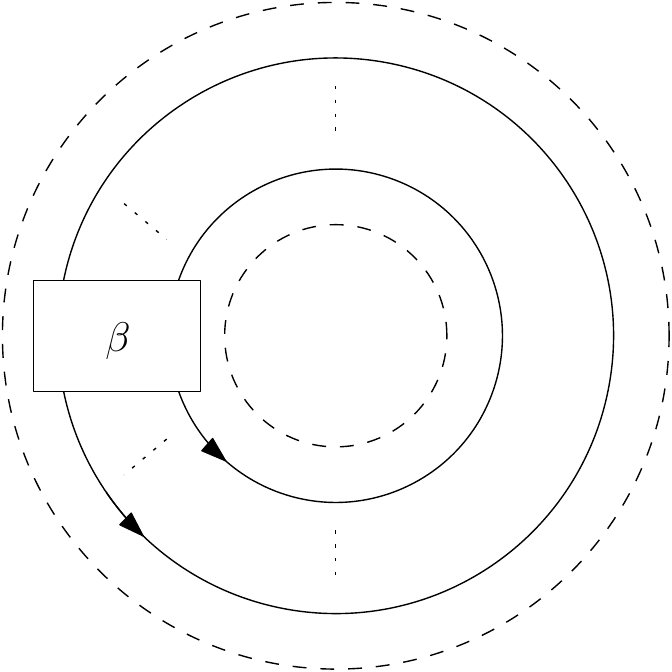}
\caption{Closure of an extended welded braid diagram.}
\label{F:closure}
\end{figure}

The \emph{closure} of an extended welded braid diagram is obtained as for usual braid diagrams (see Figure~\ref{F:closure}), with the condition that extended welded braids can be closed only when they have an even number of wen marks on each component. 

For completeness we recall that for extended welded diagrams we have two results that are  analogous to Alexander's and Markov's theorems, that we state here in the following.

\begin{prop}[{\cite[Proposition~3.3]{Damiani:ExtMarkov}}]
\label{P:Alexander}
Any extended welded link can be described as the closure of an extended welded braid diagram which is generalized Reidemeister equivalent to a welded braid diagram. 
\end{prop}

\begin{thm}[{\cite[Theorem~4.1]{Damiani:ExtMarkov}}]
\label{T:Markov}
Two extended welded braid diagrams that admit closure have equivalent closures as extended welded link diagrams if and only if they are related by a finite sequence on the following moves:
\begin{description}[before={\renewcommand\makelabel[1]{\bfseries ##1}}]
\item[$(M0)$] \label{itm:M0} isotopy of $\Rr^2$ and generalized  Reidemeister moves;
\item[$(M1)$] \label{itm:M1} conjugation in the extended welded braid group~$\WBE\nn$;
\item[$(M2)$] \label{itm:M2} a right stabilization of positive, negative or welded type, and its inverse operation.
\end{description}
\end{thm}

The $Tube$ map we briefly discussed in Subsection~\ref{SS:ext_welded_braids} can be defined also from extended welded links to ribbon torus-links, thanks to the intermediate passage through broken surfaces, and to the fact that the map is defined locally, for details see~\cite[Section~6.3]{Damiani:Journey}. On link-like objects there is no result stating that the map is an isomorphism, however we have the following result, which is a direct consequence of~\cite[Proposition~2.5]{Audoux:2016}.

\begin{prop}
\label{P:TubeSurjective}
The map $Tube$, defined on extended welded links, with values in the set of ribbon torus-links, is a well-defined surjective map.
\end{prop}

We will not linger on this construction, but we remark that the importance of this result is that it allows us to associate an extended welded link to every ribbon torus-link.

\subsection{Closed ribbon braids in \texorpdfstring{$V=B^3 \times S^1$}{}}
We introduce a particular kind of ribbon torus-links in the space $V=B^3 \times S^1$.

\begin{defn}
A torus-link $L$ in $V$ is called a \emph{closed $\nn$-ribbon braid} with $\nn \geq 1$ if $L$ meets each ball~$B^3 \times \{t\}$, for $t \in S^1$, transversely in $n$ circles.
\end{defn}

\begin{rmk}
Two closed ribbon braids in $V$ are isotopic if they are isotopic as oriented torus-links. This implies that the tubes don't necessarily stay transverse to the lamination during the isotopy.
\end{rmk}

\begin{rmk}
In general a torus-link in $V$ is not isotopic to a closed ribbon braid in~$V$. For instance a torus link lying inside a small $4$-ball in $V$ is never isotopic to a closed braid. 
\end{rmk}
 
\begin{defn}
\label{D:tubeclosure}
Given an $\nn$-ribbon braid $\beta$, its \emph{tube closure} is the ribbon torus-knot $\widehat{\beta}$ obtained by gluing a copy of the trivial ribbon braid $U$ along $\beta$, identifying the pair $(B^3 \times \{0\}, \partial_0 \beta)$ with $(B^3 \times \{1\}, \partial_1 U)$ and $(B^3 \times \{1\}, \partial_1 \beta)$ with~$(B^3 \times \{0\}, \partial_0 U)$. 
\end{defn}

On the diagrammatical side: an extended welded link diagram for $\widehat{\beta}$ in $S^1 \times I$ is obtained by closing a diagram for~$\beta$.

\section{A version of Markov's theorem in \texorpdfstring{$B^3 \times S^1$}{B3}}
\label{S:Markov}
In classical braid theory, closed braids in the solid torus are classified up to isotopy by the conjugacy classes of braids in~$\BB\nn$. We give here a classification of this kind for closed ribbon braids: their closures will be classified, up to isotopy in~$B^3 \times S^1$, by conjugacy classes of ribbon braids. The proof is inspired by the one given for the classical case in~\cite[Chapter~2]{Kassel-Turaev:2010}. In the following statement we will consider extended loop braids in their interpretation as braided annuli in the $4$-dimensional space, so we will use the terminology ``ribbon braids'' which is inherent to this approach.

\begin{thm}
\label{T:toro}
Let $\nn \geq 1$ and $\beta, \beta^\prime \in \rB\nn$ a pair of ribbon braids. The closed ribbon braids $\widehat{\beta}, \widehat{\beta^\prime}$ are isotopic in $B^3 \times S^1$ if and only if $\beta$ and $\beta^\prime$ are conjugate in~$\rB\nn$.
\end{thm}
\proof
We begin with the "if" part. Suppose first the case that $\beta$ and $\beta^\prime$ are conjugate in~$\rB\nn$. We recall that $\rB\nn$ is isomorphic to the group of extended welded braids~$\WBE\nn$. We call with the same name an element in $\rB\nn$ and a diagram for it as a representative of the corresponding class in~$\WBE\nn$. Conjugate elements of $\WBE\nn$ give rise to isotopic closed welded braid, which correspond to isotopic closed ribbon braids. This means that, since $\beta$ and $\beta^\prime$ are conjugate in $\WBE\nn$, $\beta^\prime= \alpha \beta \alpha\inv$ with $\alpha \in \WBE\nn$, and we have that~$\widehat{\alpha \beta \alpha\inv} = \widehat{\beta}$. To see this, it is enough to stack the diagrams of~$\alpha$, $\beta$ and~$\alpha\inv$, close the composed welded braid diagram, and push the upper diagram representing $\alpha$ along the parallel strands until $\alpha$ and $\alpha\inv$ are stacked one next to the other at the bottom of the diagram.

Let us now prove the converse, which is: any pair of ribbon braids with isotopic closures in $V=B^3 \times S^1$ are conjugate in~$\rB\nn$. Passing through the isomorphism between $\rB\nn$ and~$\PC\nn$, it will be enough to prove the following: any pair of ribbon braids with isotopic closures in $V=B^3 \times S^1$ have associated automorphisms of $\PC\nn$ that are conjugate. Set~$\overline{V}= B^3 \times \Rr$. Considering the cartesian product of  $(B^3, \id_{B^3})$ and the universal covering $(\Rr, p)$ of $S^1$ given by
\begin{align*}
p \colon \Rr&\longrightarrow S^1 \\
t &\longmapsto \exp(2 \pi i t)
\end{align*}
we obtain a universal covering~$(\overline{V}, \id_{B^3} \times p)$ of~$V$. Denote by $T$ the covering transformation
\begin{align*}
T \colon \overline{V}& \longrightarrow \overline{V} \\
(x, t)& \longmapsto (x, t+1)
\end{align*}
for all $x \in B^3$ and~$t \in \Rr$.
If $L$ is a closed $n$-ribbon braid in~$V$, then its preimage $\overline{L} \subset \overline{V}$ is a $2$-dimensional manifold meeting each $3$-ball~$B^3 \times \{t\}$, for $t \in \Rr$, transversely in $\nn$ disjoint pairwise unlinked circles. This implies that $\overline{L}$ consists of $\nn$ fillable components homeomorphic to~$S^1 \times \Rr$.

Being $L$ a closed ribbon braid, we can present it as a closure of a geometric ribbon braid $b \subset B^4 = B^3 \times I$ where we identify $\partial_0 B^4$ with~$\partial_1 B^4$. Then~$\overline{L}= \bigcup_{m \in \Zz} T^m (b)$,~\ie, we can see $\overline{L}$ as a tiling of an infinite number of copies of~$b$.

For~$\nn \geq 1$, let $C = (C_1, \dots, C_n)$ be a family of $\nn$ disjoint, pairwise unlinked, euclidean circles in~$\mathring{B^3}$, lying on parallel planes.
We consider a parametrization for $b$,~\ie, a family $\{\alpha_t \colon B^3 \to B^3\}_{t\in I}$ such that~$\alpha_0(C)=C$,~$\alpha_1 = \id_{B^3}$, all $\alpha_t$ fix $\partial B^3$ pointwise, and~$b=\bigcup_{t \in I} (\alpha_t(C), t)$ (see Lemma~\ref{L:NormalParametrization}).

\begin{figure}[hbtp]
\centering
\includegraphics[scale=0.4]{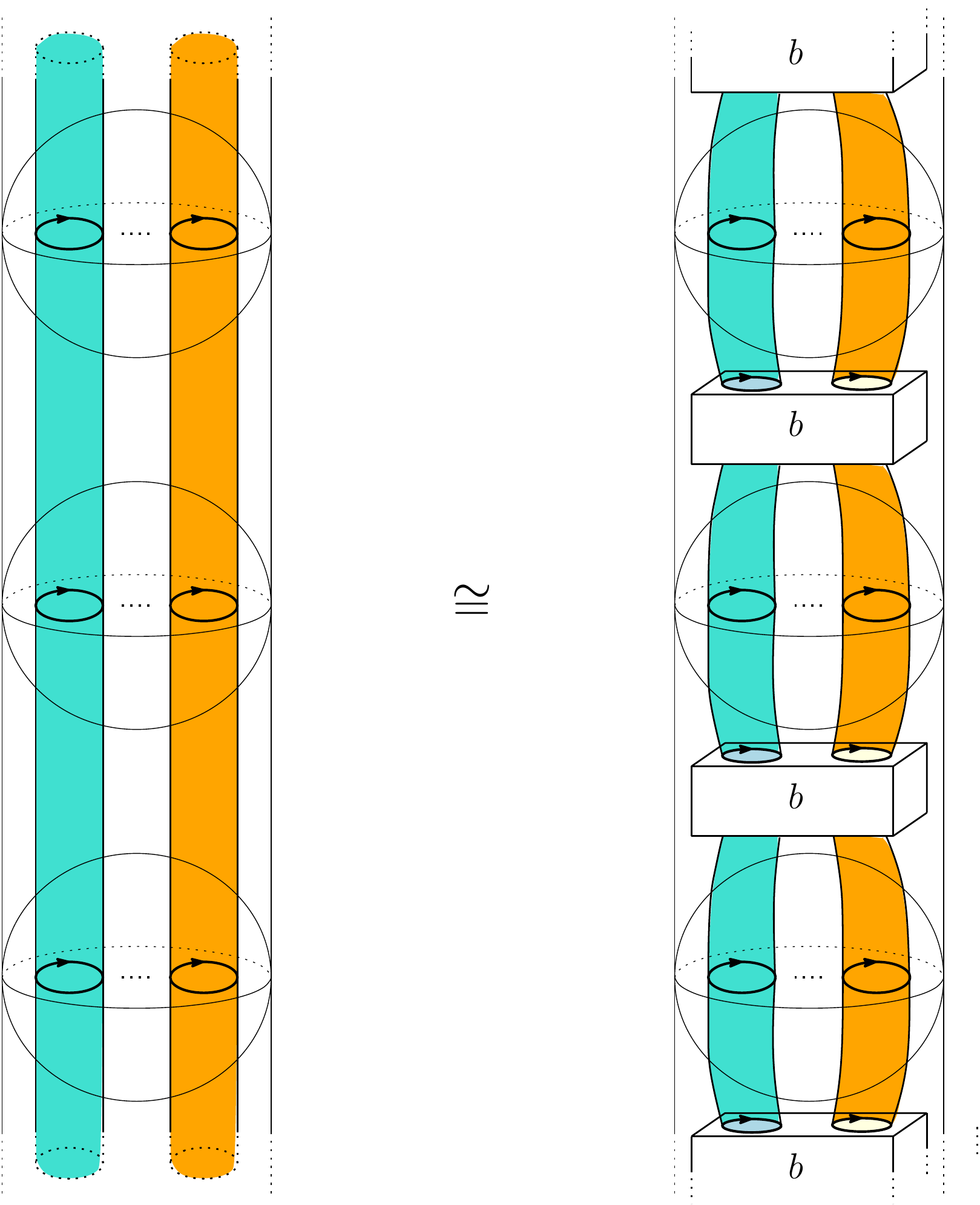}
\caption{A homeomorphism between $(B^3 \times \Rr, C \times \Rr)$ and $(B^3 \times \Rr, \overline{L})$.}
\label{F:tiling}
\end{figure}

We take the self-homeomorphism of $\overline{V}=B^3 \times \Rr$ given by 
\[
(x,t) \longmapsto (\alpha_{t-\lfloor t \rfloor} \alpha_0^{-\lfloor t \rfloor }(x), t)
\]
where~$x\in B^3$,~$t \in \Rr$, and $\lfloor t \rfloor$ is the greatest integer less than or equal to~$t$. This homeomorphism fixes $\partial\overline{V}=S^2 \times \Rr$ pointwise and sends $C\times \Rr$ onto~$\overline{L}$, see Figure~\ref{F:tiling} for an intuitive (although necessarily imprecise) idea.

The induced homeomorphism $(B^3 \setminus C)\times \Rr\cong \overline{V}  \setminus  \overline{L}$ shows that $B^3 \setminus C = (B^3 \setminus C)\times \{0\} \subset \overline{V}  \setminus  \overline{L}$ is a deformation retract of~$\overline{V}  \setminus  \overline{L}$. Pick a point $d \in \partial_0 B^4 = B^3$ and set~$\overline{d}=(d,0)\in \overline{V}$; them the inclusion homomorphism
\[
i \colon \pi_1(B^3 \setminus C, d) \longrightarrow \pi_1(\overline{V}  \setminus  \overline{L}, \overline{d})
\]
is an isomorphism.

By definition the image of $d$ by the covering transformation $T$ is~$T(d)=(d, 1)$; the covering transformation $T$ restricted to $\overline{V}  \setminus  \overline{L}$ induces an isomorphism~$\pi_1(\overline{V}  \setminus  \overline{L},\overline{d}) \to \pi_1(\overline{V}  \setminus  \overline{L},T(\overline{d}))$.
Let $T_\ast$ be the composition of this isomorphism with the isomorphism $\pi_1(\overline{V}  \setminus  \overline{L},T(\overline{d})) \to \pi_1(\overline{V}  \setminus  \overline{L},\overline{d})$ obtained by conjugating the loops by the path~$d \times [0,1] \subset \partial B^3 \times \Rr\subset \overline{V}  \setminus  \overline{L} $. Then $T_\ast$ is an automorphism of $\pi_1(\overline{V}  \setminus  \overline{L},\overline{d})$. Therefore the following diagram commutes:
\[
\begin{CD}
\pi_1(B^3 \setminus C, d) @>i>> \pi_1(\overline{V}  \setminus  \overline{L},\overline{d}) \\
@V{\tilde{\beta}}VV  @VV{T_\ast}V @. \\
\pi_1(B^3 \setminus C, d) @>i>> \pi_1(\overline{V}  \setminus  \overline{L},\overline{d})\\
\end{CD}
\]
where $\tilde{\beta}$ is the automorphism induced by the restriction of $\alpha_0$ to~$B^3 \setminus C$.
The isomorphism between $\rB\nn$ and $\MCG {B^3}{C^\ast}$ allows us to send the ribbon braid~$\beta$, represented by~$b$, to the isotopy class of~$\alpha_0$.

Identifying $\pi_1(B^3 \setminus C, d)$ with the free group $\F\nn$ with generators $x_1, x_2, \dots, x_n$, we conclude that the automorphism~$\tilde{\beta}$ is equal to~$\nu(b)$, where $\nu \colon \rB\nn \to \PC\nn$ is the isomorphism between the group of ribbon braids $\rB\nn$ and~$\PC\nn$, the subgroup of~$\Aut(\F\nn)$ generated by the automorphisms of the form $\alpha \colon x_\ii \mapsto w_\ii\inv x^{\pm 1}_{\pi({\ii})} w_\ii$ where $\pi$ is a permutation and $w_\ii$ is a word in~$\F\nn$. Then it is the automorphism of $\F\nn$ corresponding to~$\beta$, the ribbon braid represented by~$b$. Thus~$i\inv T_\ast i = \tilde{\beta}$.

Suppose now that $\beta, \beta^\prime \in \rB\nn$ are two ribbon braids with isotopic closures in~$V$, and that $b$ and $b^\prime \subset B^4 =B^3\times I$ are two geometric ribbon braids that represent them. Let $L$ and $L^\prime \subset V = B^3 \times S^1$ be their respective closures. 

Then there is a homeomorphism $g \colon V \to V$ such that $g$ maps $L$ onto $L^\prime$, preserving their canonical orientation along the annuli, but possibly reversing the orientation of the circles at some instant (for example when Reideiester moves of type I occur). Note that a Reidemeister move of type I is isotopic to the composition of two wens~\cite[Corollary 3.3]{Audoux:2016}, so globally the orientation of the circles at the starting and ending time parameter is preserved). In fact the orientation of the ambient $V$ is preserved by~$g$, but when considering a section $B^3 \times \{t\}$ the orientation of the circles can be concordant or not concordant with the one induced by~$V$.
In addition the restriction of $g$ to $\partial V$ is isotopic to the identity~$\id_V$. This fact, plus the isomorphism of the map induced by the inclusion $\pi_1(\partial V)=\pi_1(S^2 \times S^1) \to \pi_1(V)=\pi_1(B^3 \times S^1) \cong \Zz$ implies that $g$ induces an identity automorphism of~$\pi_1(V)$. Therefore $g$ lifts to a homeomorphism $\overline{g} \colon \overline{V} \to \overline{V}$ such that $\overline{g}$ is isotopic to the identity on $\partial \overline{V}$, $\overline{g}T=T \overline{g}$, and~$\overline{g}(\overline{L})=\overline{L^\prime}$.

Hence $\overline{g}$ induces an isomorphism 
\[
\overline{g}_\ast \colon \pi_1(\overline{V} \setminus \overline{L}, \overline{d}) \longrightarrow\pi_1(\overline{V} \setminus \overline{L^\prime}, \overline{d})
\]
commuting with~$T_\ast$. The following diagram commutes:
\[
\begin{CD}
\pi_1(B^3 \setminus C, d) @>i>> \pi_1(\overline{V}  \setminus  \overline{L},\overline{d}) \\
@V{\varphi}VV  @VV{\overline{g}_\ast}V @. \\
\pi_1(B^3 \setminus C, d) @>i^\prime>> \pi_1(\overline{V}  \setminus  \overline{L^\prime},\overline{d}).\\
\end{CD}
\]

Consider the automorphism $\varphi = (i^\prime)\inv \overline{g}_\ast i$ of~$\F\nn=\pi_1(B^3 \setminus C, d)$, where:
\begin{align*}
&i \colon \pi_1 (B^3 \setminus C, d) \longrightarrow \pi_1(\overline{V} \setminus \overline{L}, \overline{d}) \ \mbox{and }\\
&i^\prime \colon \pi_1 (B^3 \setminus C, d) \longrightarrow \pi_1(\overline{V} \setminus \overline{L^\prime}, \overline{d})
\end{align*}
are the inclusion isomorphisms.

Applying the same arguments to~$\beta^\prime$, we have $\tilde{\beta^\prime} = (i^\prime)\inv T_\ast i^\prime$, and from the preceding commutative diagram we have: 
\[
\varphi \tilde{\beta} \varphi\inv= \big( ({i^\prime})\inv  \overline{g}_\ast i \big) \ \big( i\inv T_\ast i \big) \ \big( i\inv {\overline{g}_\ast}\inv i^\prime \big) = (i^\prime)\inv T_\ast i^\prime =  \tilde{\beta^\prime}
\]

We claim that $\varphi$ is an element of the subgroup of $\Aut(\F\nn)$ consisting of all automorphisms of the form $x_i \mapsto q_i x_{j(i)}^{\pm1}q_i\inv$, where~$i=1, \dots, n$, $j(i)$ is some permutation of the numbers~$1, \dots, n$, and $q_i$ a word in~$x_1, \dots, x_n$.
Then the isomorphism between this subgroup and~$\rB\nn$ implies that $\beta$ and $\beta^\prime$ are conjugate in~$\rB\nn$.

We prove this claim. The conjugacy classes of the generators $x_1, x_2, \dots, x_n$ in $\F\nn = \pi_1(B^3 \setminus C, d)$ are represented by loops encircling the circles~$C_i$. The inclusion $B^3 \setminus C=(B^3 \setminus C)\times \{0\} \subset \overline{V} \setminus \overline{L}$ maps these loops to some loops in $\overline{V} \setminus \overline{L}$ encircling at each parameter $t$ the rings that form the components of~$\overline{L}$. 
The homeomorphism $\overline{g}\colon \overline{V} \to \overline{V}$ transforms these loops into loops in $\overline{V} \setminus \overline{L^\prime}$ encircling the components of~$\overline{L^\prime}$.
The latter represent the conjugacy classes of the images of $x_1, \dots, x_n$ under the inclusion~$B^3 \setminus C=(B^3 \setminus C) \times \{0\} \subset \overline{V} \setminus \overline{L^\prime}$.

The automorphism $\varphi$ transforms the conjugacy classes of $x_1, \dots, x_n$ into themselves, up to permutation and orientation changes. This verifies the condition. The possible orientation changes are due to the fact that the isotopy of closed braid is not monotone with respect to the time parameter as ribbon braid isotopy is, thus Reidemeister moves of type I can occur. 
\endproof

When one ribbon braid is a conjugate of another ribbon braid, we can describe the form of the conjugating element.

\begin{lem}
Let $\nn \geq 1$ and $\beta, \beta^\prime \in \rB\nn$ a pair of ribbon braids. They are conjugates in $\rB\nn$ if and only if $\beta^\prime = \pi_{\tau} \alpha \beta \alpha\inv \pi_{\tau}\inv$, where $\pi_\tau$ is composed only by wens and~$\alpha$ does not contain any wen. Speaking in terms of presented group, $\pi_\tau$  is represented by a word in the $\tau_\ii$ generators of presentation~\eqref{E:Rpresentation}.
\end{lem}
\proof
Take $\beta$ and $\beta^\prime$ in $\rB\nn$ conjugate by another element in~$\rB\nn$. Then there exists an element $\gamma$ in $\rB\nn$ such that~$\beta = \gamma \beta^\prime \gamma\inv$. Consider $\gamma$ as an element of the configuration space of $\nn$ circles~$\R\nn$. We can use relations from presentation \eqref{E:Rpresentation} to push to the right of the word the  generators~$\tau_\ii$, to obtain an equivalent element~$\gamma\prime = \pi_\tau \alpha$, where $ \pi_\tau$ is a word in the $\tau_\ii$s and $\alpha$ only contains generators  $\sig\ii$ and~$\rr\ii$. This means that $\alpha$ is in fact an element that can be written without $\tau$ generators. Finally, when considering $\gamma\inv$ for the conjugacy, we remark that $\pi_\tau\inv$ is just the mirror image word of~$\pi_\tau$. 
\endproof

\section{Ideas for further developements}
\label{S:future}
To extend the result in~$\Rr^4$  we shall prove the invariance of isotopy classes of closed ribbon braids under the operation known as stabilisation.
The approches used for usual knotted objects, for instance those of \cite{Birman:Libro} and~\cite{Traczyk:Markov}, rely on the bijection  given by Reidemeister theorem between knots and knot diagrams up to Reidemeister moves. We do not have such a result for ribbon torus-links. In fact, as Proposition~\ref{P:TubeSurjective} points out, the injectivity of the map $Tube$ between extended welded links and ribbon torus-links is an open question. When applied to welded links (not extended), we know that the $Tube$ map is not injective: for instance, it is invariant under the \emph{horizontal mirror image} on welded diagrams (\cite[Proposition~3.3]{Ichimori-Kanenobu:2012}, see also \cite{Winter:2009, Satoh:2000}), while welded links in general are not equivalent to their horizontal mirror images. 
However, extended welded links \emph{are} equivalent to their horizontal mirror image~(\cite[Proposition~5.1]{Damiani:ExtMarkov}). This fact suggest they could be good candidates to be in bijection with ribbon torus-links. Of course, other obstructions to injectivity may exist, so the relation between extended welded links and ribbon torus-links shall be investigated. It is worth noticing that an alternative approach to solve the problem of establishing a bijection between welded diagrams and ribbon torus-links has been suggested by Kawauchi in~\cite[Problem, Section~2]{Kawauchi:Chord16moves}.

\section*{Acknowledgements}
During the writing of this paper, the author was supported by a JSPS Postdoctral Fellowship For Foreign Researchers and by JSPS KAKENHI Grant Number 16F16793. The author thanks Emmanuel Wagner for valuable conversations.


\bibliography{KIH_Markov.bib}{}
\bibliographystyle{alpha}

\end{document}

%% file: Macro.tex
\newcommand\Nn{\mathbb{N}}
\newcommand\Zz{\mathbb{Z}}

\newcommand\Rr{\mathbb{R}}

\DeclareMathOperator{\Aut}{Aut}	

\newcommand\BB[1]{B_{#1}}	
\newcommand\WB[1]{WB_{#1}}	
\newcommand\WBE[1]{WB_{#1}^{ext}}	
\newcommand\PC[1]{PC_{#1}}	
\newcommand\PLBE[1]{PLB_{#1}^{ext}}	
\newcommand\LBE[1]{LB_{#1}^{ext}}	
\newcommand\LB[1]{LB_{#1}}	
\newcommand\rB[1]{rB_{#1}}	

\newcommand \ep{\varepsilon}

\newcommand \F[1]{F_{#1}}

\DeclareMathOperator{\Homeo}{Homeo}	
\DeclareMathOperator{\Diffeo}{Diffeo}	
\DeclareMathOperator{\PHomeo}{PHomeo}	

\newcommand\id{\mathrm{id}}
\newcommand\ie{{\textit{i.e.}}}
\newcommand\ii{i}
\newcommand\inv{^{-1}}

\newcommand\jj{j}

\newcommand\MCG[2]{\mathrm{MCG}({#1}, {#2})} 	
\newcommand\PMCG[2]{\mathrm{PMCG}({#1}, {#2})} 	

\newcommand\nn{n}
\newcommand\nno{{n-1}}

\makeatletter
\newcommand{\oset}[2]{%
  {\mathop{#2}\limits^{\vbox to -.5\ex@{\kern-\tw@\ex@
   \hbox{\scriptsize #1}\vss}}}}
\makeatother

\newcommand\R[1]{R_{#1}}	
\newcommand\Ri[1]{\mathcal{R}_{#1}}	

\newcommand\rr[1]{\rho_{\hspace{-0.3ex}#1}^{\null}}	
\newcommand\rrq[2]{\rho_{\hspace{-0.3ex}#1}^{#2}}	

\newcommand\sig[1]{\sigma_{\hspace{-0.3ex}#1}^{\null}} 
\newcommand\sigg[2]{\sigma_{\hspace{-0.3ex}#1}^{#2}}	
\newcommand\siginv[1]{\sigg{#1}{-1}}	
